\documentclass[final,leqno,letterpaper]{myetna}
\usepackage{amsmath}
\usepackage{amssymb}
\usepackage{xspace}
\usepackage{xcolor}

\newtheorem{assumption}[theorem]{Assumption}
\newtheorem{algorithm}[theorem]{Algorithm}

\title{Convergence results and low order rates for nonlinear Tikhonov regularization with oversmoothing penalty term}

\author{Bernd Hofmann\footnotemark[2]
		 \and Robert Plato\footnotemark[3]}

\newcommand{\para}{\alpha}
\newcommand{\parac}{\beta}
\newcommand{\pardel}{\para_*}
\newcommand{\pardelb}{\beta_*}
\newcommand{\ix}{X}
\newcommand{\yps}{Y}

\newcommand{\fdelta}{f^\delta}
\newcommand{\reza}{ \mathbb{R} }
\newcommand{\Landau}{\mathcal{O}}
\newcommand{\refeq}[1]{(\ref{eq:#1})}

\newcommand{\normqua}[1]{\norm{#1}^2}

\newcommand{\norm}[1]{\Vert \hspace{0.4mm} #1 \hspace{0.4mm} \Vert}

\newenvironment{myenumerate}{%
\begin{list}{(\alph{enumcount})}
{\setcounter{enumcount}{1}\usecounter{enumcount}
\setlength{\topsep}{1mm}
\setlength{\itemsep}{0mm}
\setlength{\labelwidth}{0mm}
\setlength{\labelsep}{1mm}
\setlength{\itemindent}{1mm}
\setlength{\leftmargin}{0mm}
}}{\end{list}}

\newcommand{\DB}{\Dset(B)}
\newcommand{\DF}{\Dset(\F)}

\newcommand{\Dset}{\mathcal{D}}

\newcommand{\myu}{u}
\newcommand{\myv}{v}

\newcommand{\ust}{u^\dagger}

\newcommand{\fst}{f^\dagger}
\newcommand{\fdel}{f^\delta}
\newcommand{\upardeldel}{u_{\pardel}^\delta}

\newcommand{\tikreg}{Tikhonov regularization\xspace}
\newcommand{\tifu}{Tikhonov functional\xspace}

\newcommand{\tinybullet}{{\tiny \raisebox{0.6mm}{$ \bullet $}}}
\newenvironment{mylist}{%
\begin{list}{\tinybullet}
{\setlength{\topsep}{0.2cm}
\setlength{\itemsep}{0mm}
\setlength{\labelwidth}{0mm}
\setlength{\labelsep}{3mm}
\setlength{\itemindent}{3mm}
\setlength{\leftmargin}{0mm}
}}{\end{list}}

\newenvironment{mylist_indent}{%
\begin{list}{\tinybullet}
{\setlength{\topsep}{0.2cm}
\setlength{\itemsep}{0mm}
\setlength{\labelwidth}{2mm}
\setlength{\labelsep}{3mm}
\setlength{\itemindent}{0mm}
\setlength{\leftmargin}{5mm}
}}{\end{list}}

\newcommand{\kla}[1]{(#1)}
\newcommand{\mfrac}[2]{\dfrac{\mbox{\footnotesize \raisebox{-0.5mm}{$#1$}}}%
{\mbox{\footnotesize \raisebox{0.8mm}{$#2$}}}}
\newcommand{\proofend}{\qquad \endproof}

\newcommand{\Landauno}[1]{\Landau\kla{#1}}
\newcommand{\as}{\quad \text{as} \ \ }
\newcommand{\assh}{\ \  \text{as} \ }

\newcommand{\R}{\mathcal{R}}
\newcommand{\cf}{cf.\mbox{}\xspace}
\newcommand{\ie}{i.e.,\xspace}

\newcommand{\bn}{\bigskip \noindent}
\newcommand{\defeq}{:=}
\newcommand{\for}{\quad \text{for} \ \ }

\newcommand{\rhs}{right-hand side\xspace}

\newcommand{\parab}{\beta}

\newenvironment{myenumerate_indent}{%
\begin{list}{(\alph{enumcount})}
{\setcounter{enumcount}{1}\usecounter{enumcount}
\setlength{\topsep}{1mm}
\setlength{\itemsep}{0mm}
\setlength{\labelwidth}{5mm}
\setlength{\labelsep}{2mm}
\setlength{\itemindent}{-0mm}
\setlength{\leftmargin}{7mm}
}}{\end{list}}

\newcommand{\inset}[1]{\{ \, #1 \, \}}

\newcommand{\myast}{\cdot}

\newcommand{\modul}[1]{\vert \hspace{0.3mm} #1 \hspace{0.3mm} \vert}

\newcommand{\ualpaux}[1][\para]{\widehat{u}_{#1}}
\newcommand{\Hop}{G}

\newcommand{\fab}{2a+2}
\newcommand{\norma}[1]{\norm{#1}_{-a}}

\newcommand{\normone}[1]{\norm{#1}_{1}}
\newcommand{\normonequa}[1]{\norm{#1}_{1}^2}
\newcommand{\normtau}[2][\tau]{\norm{#2}_{#1}}
\newcommand{\normyps}[1]{\norm{#1}}
\newcommand{\normypsqua}[1]{\norm{#1}^2}
\newcommand{\normix}[1]{\norm{#1}}
\newcommand{\landau}{o}

\newcommand{\landauno}[1]{\landau\kla{#1}}
\newcommand{\lfrac}[2]{#1/#2}
\newcommand{\mytheta}{\theta}
\newcommand{\ca}{c_a}
\newcommand{\cb}{C_a}
\newcommand{\upardel}[1][\para]{{u}_{#1}^\delta}
\newcommand{\upardelb}{\upardel[\parab]}
\newcommand{\updd}{\upardel[\pardel]}
\newcommand{\upddb}{\upardel[\pardelb]}
\newcommand{\udel}{\updd}
\newcommand{\tikfun}{T}
\newcommand{\pad}[2][\para]{\tikfun_{#1}^\delta(#2)}
\newcommand{\padpur}[1][\para]{\tikfun_{#1}^\delta}
\newcommand{\Jdelqua}[1]{\normqua{F{#1} - \fdel}}
\newcommand{\Jdel}[1]{\norm{F{#1} - \fdel}}
\newcommand{\J}[1]{\norm{F{#1} - \fst}}
\newcommand{\myomega}[1]{\normonequa{#1}}
\newcommand{\F}[1]{\Fpur#1}
\newcommand{\Fpur}{F}
\newcommand{\myb}{b}
\newcommand{\myc}{c}
\newcommand{\uadno}[1][k]{u_{\alpha^{(#1)}}^\delta}
\newcommand{\parainit}{\para^{(0)}}

\newcommand{\ubar}{\overline{u}}

\newcommand{\mainassump}{Let Assumption \ref{th:main_assump} be satisfied.\xspace}
\newenvironment{myenumerate_roman}{%
\begin{list}{(\roman{enumcountroman})}
{\setcounter{enumcountroman}{1}\usecounter{enumcountroman}
\setlength{\topsep}{0.2cm}
\setlength{\itemsep}{0mm}
\setlength{\labelwidth}{0mm}
\setlength{\labelsep}{3mm}
\setlength{\itemindent}{3mm}
\setlength{\leftmargin}{0mm}
}}{\end{list}}

\newcommand{\wrt}{with respect to\xspace}

\newcommand{\myH}{\Delta}
\newcommand{\myHp}{\Delta^\prime}
\newcommand{\nondecreasing}{non-decreasing\xspace}
\newcommand{\monotonically}{}
\newcommand{\strictly}{strictly\xspace}
\newcommand{\nonincreasing}{non-increasing\xspace}
\newcommand{\rp}[1]{\textcolor{black}{#1}}
\newcommand{\rpc}[1]{#1}
\newcommand{\rptwo}[1]{\textcolor{magenta}{#1}}

\newcommand{\bh}[1]{\textcolor{black}{#1}}
\newcommand{\elltwo}{\ell_2}
\newcommand{\normb}[1]{\norm{#1}_{\elltwo}}
\newcommand{\normbqua}[1]{\normb{#1}^2}
\newcommand{\N}{N}
\newcommand{\aninv}{n}
\newcommand{\an}{\tfrac{1}{n}}
\newcommand{\mykap}{\kappa}
\newcommand{\myq}{2.31}

\makeatletter
\def\proof{\par{\textit{Proof. }}\ignorespaces\global\@endproofsetfalse}
\makeatother

\shorttitle{TIKHONOV REGULARIZATION WITH OVERSMOOTHING PENALTY TERM}
\shortauthor{B.~HOFMANN AND R.~PLATO}

\begin{document}
\maketitle

\renewcommand{\thefootnote}{\fnsymbol{footnote}}

\footnotetext[2]{Faculty of Mathematics, Chemnitz University of Technology, 09107 Chemnitz, Germany.}
\footnotetext[3]{Department of Mathematics, University of Siegen,
Walter-Flex-Str.~3, 57068 Siegen, Germany}

\newcounter{enumcount}
\newcounter{enumcountroman}
\renewcommand{\theenumcount}{(\alph{enumcount})}
\bibliographystyle{plain}
\begin{abstract}
\rp{For the Tikhonov regularization of ill-posed nonlinear operator equations,
convergence is studied in a Hilbert scale setting.}
We include the case of oversmoothing penalty terms, which means that the exact solution does not belong to the domain
of definition of the considered penalty functional. In this case, we try to close a gap in the present theory,
where H\"older-type convergence rates results have been proven under corresponding source conditions, but assertions on norm convergence of
regularized solutions without source conditions are completely missing.
A result of the present work is to provide sufficient conditions for convergence under a priori and a posteriori regularization
parameter choice strategies, without any additional smoothness assumption on the solution. The obtained error estimates moreover allow us to prove
low order convergence rates under associated (for example logarithmic) source conditions.
\rp{Some numerical illustrations are also given.}
\end{abstract}


\section{Introduction}
\label{intro}

The subject of this paper are nonlinear operator equations of the form
\begin{equation} \label{eq:opeq}
 \F u = \fst \,,
\end{equation}
where $ \F: \ix \supset \DF \to \yps $ is a nonlinear operator
between infinite-dimensional Hilbert spaces $ \ix $ and $ \yps $ with norms $\|\cdot\|$.
We suppose that the \rhs $ \fst \in Y $ is approximately given as $ \fdelta \in \yps$ satisfying the deterministic noise model
\begin{equation} \label{eq:noise}
 \normyps{ \fdelta-\fst } \le \delta,
\end{equation}
with the noise level $\delta \ge 0$.
Throughout the paper, it is assumed that the considered equation~\eqref{eq:opeq} has a solution $ \ust \in \DF $ and is (at least locally at $ \ust$) ill-posed (cf.~\cite{HofPla18}).

For finding stable approximations to the solution  $ \ust \in \DF $ of equation \eqref{eq:opeq}, we consider the \tikreg, where the regularized solutions are minimizers of the extremal problem
\begin{equation} \label{eq:TR}
 \pad{\myu} \defeq \Jdelqua{\myu} + \para \myomega{\myu-\ubar} \to \min \quad \mbox{subject to} \quad  \myu \in \DF,
\end{equation}
with a regularization parameter $ \para > 0 $.
In this context, $ \normone{\cdot} $ is assumed to be a norm of a densely defined subspace $X_1$ of $X$, which is stronger than the original norm $ \normix{\cdot} $ in $X$. Throughout this paper, we suppose that
the initial guess $\ubar$ occurring in the penalty term of $\pad{\myu}$ satisfies the condition
\begin{equation} \label{eq:guess}
\ubar \in \Dset:=\DF \cap X_1.
\end{equation}

Precisely, we define the stronger norm $\|\cdot\|_1$ by a generator $ B: \ix \supset \DB \to \ix $, which is a self-adjoint and positive definite unbounded linear operator with dense domain $ \DB $, i.e.~we have for some constant $ m > 0 $
\begin{align}
\normix{B\myu} \ge m \normix{ \myu } \for \myu \in \DB.
\label{eq:b-inverse-stable}
\end{align}
This allows us to introduce norms
\begin{align}
\label{eq:normtau}
\normtau{ \myu } \defeq \norm{B^\tau\myu}, \quad \myu \in \ix_\tau
\qquad (\tau \in \reza),
\end{align}
where $ \ix_\tau := \Dset(B^\tau) $ for $ \tau > 0 $, and
$ \ix_\tau := \ix $ for $ \tau \le 0 $.
The fractional powers are defined by means of the resolution of the identity generated by the inverse operator $ B^{-1} $,
see, e.g., \cite[Section 2.3]{EHN96}.
Note that the system of spaces $ (\ix_\tau)_{\tau \in \reza} $, equipped with the respective norms, is strongly related to the Hilbert scale generated by the operator $ B $. However, for $ \tau < 0 $, topological completion of the spaces $ \ix_\tau = \ix $ \wrt the norm $ \normtau{\cdot} $
is not needed in our setting, and thus is omitted.

In the present work, we discuss the nonlinear \tikreg \eqref{eq:TR} in particular with an oversmoothing penalty term, where we have $ \ust \not \in \ix_1=\DB$, or in other words $\|\ust\|_1=+\infty.$
This continues studies started in papers \cite{Hofmann_Mathe[19],Hofmann_Mathe[18]} and \cite{GHH19}, where
convergence rates and numerical case studies are provided for a priori and a posteriori parameter choices, respectively, under certain smoothness assumptions on $ \ust $ and structural conditions on $\F$.
Under the same structural conditions, which are also similar to those in the corresponding seminal paper for linear operator equations by Natterer \cite{Natterer[84]}, we present as the novelty of this paper convergence results  based on the Banach--Steinhaus theorem without needing any smoothness assumptions. The error estimates derived in the context of convergence assertions moreover allow us to prove
low order convergence rates under associated (for example logarithmic) source conditions.

The outline of the remainder is as follows: in Section~\ref{preparations},
we introduce Hilbert scales and formulate the basic assumptions, and in addition we establish well-posedness of \tikreg used in our setting.
Then in Section~\ref{auxels}, we introduce auxiliary elements needed for the proof of the convergence results, and in addition we provide first error estimates for \tikreg which are based on those auxiliary elements and which are needed for the subsequent convergence proofs.
The regularizing properties of an a priori parameter choice as well as a discrepancy principle are considered  in Section~\ref{parameter_choices}. The suggested discrepancy principle is considered in a form that is suitable for misfit functionals which may depend discontinuously on the regularization parameter $ \para>0 $.
As a byproduct of derived error estimates, we can prove low order convergence rates in Section~\ref{sec:lowrates}.
\rp{We conclude this paper by presenting results of some numerical experiments.}
\section{Prerequisites and assumptions}
\label{preparations}
\subsection{Main assumptions}
In the following assumption, we briefly summarize the structural properties of the operator $F$ \rp{and} of its domain $\DF$, in particular with respect to the the solution $\ust$ of equation
\eqref{eq:TR}. For examples of nonlinear inverse problems, which satisfy these assumptions (or at least substantial parts of it), we refer to \cite{Egger_Hofmann[18],GHH19} and to the appendices of the papers \cite{Hofmann_Mathe[18],WerHof20}.
\begin{assumption}
\label{th:main_assump}
\begin{itemize}
\item[(a)]
The operator $ \F: \ix \supset \DF \to \yps $ is sequentially continuous
on $ \DF $ \wrt the weak topologies of the Hilbert spaces $ \ix $ and $ \yps $.

\item[(b)]
The domain of definition $ \DF \subset \ix $ is a closed and convex subset of $X$.
\item[(c)]
Let $\Dset \defeq \DF \cap \DB=\DF \cap X_1 $ be a non-empty set.

\item[(d)]
Let the solution $ \ust \in \DF $ to equation \eqref{eq:opeq} with \rhs $\fst$ be
an interior point of the domain $ \DF $.

\item[(e)]
Let the data $ \fdelta \in \yps$ satisfy the noise model \eqref{eq:noise}, and let the
initial guess $\ubar$ satisfy \eqref{eq:guess}.

\item[(f)] Let $ a > 0 $, and let there exist finite constants $ 0 < \ca \le \cb $ such that
the inequality chain
\begin{align}
\label{eq:normequiv}
\ca\norma{\myu - \ust} \le \normyps{\F{\myu} - \fst} \le
\cb\norma{\myu - \ust}
\end{align}
holds true for all $\myu \in \Dset$.
\end{itemize}
\end{assumption}
\begin{remark} \label{rem:rem1}
From item (f) (left-hand inequality) of Assumption~\ref{th:main_assump}, we have for $\ust \in X_1$ that $\ust$ is the uniquely determined solution to equation \eqref{eq:opeq} in the set $\Dset$.
For $\ust \notin X_1$, there is no solution at all to \eqref{eq:opeq} in $\Dset$. But in both cases, alternative solutions $u^* \notin X_1$ with $u^* \in \DF$ and $F u^*=\fst$ cannot be excluded.
\end{remark}

\subsection{Properties of regularized solutions of the Tikhonov regularization}
\rp{For $\alpha>0$, minimizers of the Tikhonov functional $\padpur $
\rp{exist (cf.~Proposition \ref{th:tifu-well-posed} below) and}
are denoted by $ \upardel$, i.e.~we have $\pad{\upardel} = \min \nolimits_{\myu \in \DF} \ \pad{u}$.
Evidently, by definition of the penalty term,
$ \upardel \in \Dset$ holds.}

\begin{example}
\label{th:linear} Let in this example $ \Fpur = A: \ix \to \yps $ with $\DF=\ix$ be a bounded linear operator with non-closed range $\mathcal{R}(A)$, and for simplicity let $ \ubar = 0 $. In this setting,
Tikhonov regularized solutions $\upardel$ solve the linear operator equation
\begin{equation} \label{eq:auxFormel}
(A^* A + \para B^{-2}) \,\upardel = A^* \fdelta.
\end{equation}
In the special situation of an injective operator $ A $ and of a scale generator $ B = (A^*A)^{-q/2} $ with $ q > 0 $, this gives
$$ (A^* A + \para (A^*A)^q)\, \upardel = A^* \fdelta ,$$
and Assumption \ref{th:main_assump}
is satisfied with $ a = 1/q $ then.
The oversmoothing case $ \ust \not \in \ix_1 $ here means
$ \ust \not \in \mathcal{R}((A^*A)^{q/2}) $.
This situation is discussed in the analysis of fractional \tikreg, and we refer for example to
\cite{Bianchi_Buccini_Donatelli_SerraCapizzano[15],Hochstenbach_Reichel[11],Louis[89]}. In Natterer's paper \cite{Natterer[84]}, the analog
\begin{equation} \label{eq:Nata}
\ca\norma{\myu} \le \normyps{A \myu} \le \cb\norma{\myu} \qquad \mbox{for all} \quad \myu \in X
\end{equation}
to the inequality chain \eqref{eq:normequiv} is the basis for error estimates and convergence rates results for linear operator equations. The constant $a>0$  characterizes here the degree of ill-posedness of the problem.
\bh{We mention that Neubauer has discussed in \cite{Neubau92} the consequences of the two-sided condition
\begin{equation} \label{eq:Neuba}
\ca\norma{\myu} \le \normyps{F^\prime(\ust){\myu}} \le \cb\norma{\myu} \qquad \mbox{for all} \quad \myu \in X,
\end{equation}
which is an extension of (\ref{eq:Nata}) to the nonlinear case and closely connected to (\ref{eq:normequiv}), if the forward operator $F$ possesses a Fr\'echet derivative $F^\prime(\ust)$ at $\ust$. For
more details, see also \cite[Sect.~4]{Egger_Hofmann[18]}}.
\end{example}

The extremal problem \eqref{eq:TR} for finding regularized solutions is well-posed with respect to existence of minimizers and their stability, in a sense specified in the following proposition. This follows by standard results
from regularization theory (cf., e.g., \cite[Chapter 2.6]{Tikhonov_Leonov_Yagola[98]}, \cite{Vainikko[88],Vainikko[91.1]} and \cite[Section~4.1.1]{Schusterbuch12}). So we give a sketch of proof only.
\begin{proposition}
\label{th:tifu-well-posed}
\mainassump
\begin{myenumerate_indent}
\item
There exists, for all $\alpha>0$, a minimizer $ \upardel$ of the \tifu $\padpur$ in the set $\Dset $.
\item
Each minimizing sequence of $\padpur$ over $\Dset$ has a subsequence that converges strongly in $ \ix_1 $ to a minimizer $ \upardel \in \Dset $ of the Tikhonov functional.
\item For every $\alpha>0$, the regularized solutions  $ \upardel$ are stable \rptwo{in $ X_1 $} with respect to small perturbations in the data $\fdelta$.
\end{myenumerate_indent}
\end{proposition}
\proof The basic ingredients needed for the proof are as follows:
\begin{mylist_indent}
\item
The operator $ \F $, when considered as $ \F: \ix_1 \supset \Dset \to \yps $, is sequentially continuous \wrt the weak topologies on $ \ix_1 $ and $ \yps $.
This implies that the misfit functional
$\myu \in \Dset \mapsto \Jdel{u} \in \reza$
is sequentially continuous \wrt the weak topology on $ \ix_1 $.

\item The set $ \Dset $ is weakly closed in $ \ix_1 $.

\item
The stabilizing functional $ \normone{\cdot - \ubar}^2 $ is sequentially weakly lower continuous
on $ \ix_1 $.
\end{mylist_indent}
The statement in the second item follows from the two facts that (i)
\rp{the embedding operator $ \ix_1 \hookrightarrow \ix $ is continuous},
and that (ii) each closed convex subset of a Hilbert space is weakly closed.

From these ingredients, it follows that each minimizing sequence $ (u_n) \subset \Dset $ of the \tifu has a subsequence which converges weakly in $ \ix_1 $ to a minimizer $ \upardel $, and the corresponding subsequence of $ (\normone{u_n - \ubar}) $ converges to $ \normone{\upardel - \ubar} $.
\proofend

\begin{remark}
\bh{We note that the minimizer of the \tifu may be non-unique,
because $T_\alpha^\delta$ can be, for nonlinear forward operators $F$,
a non-convex functional as a consequence of a non-convex misfit term
$\|F u - f^\delta\|^2$. If, for example, $Fu:=u \star u$ represents the autoconvolution
operator in $X=L^2(0,1)$ (cf., e.g.,~\cite{BueHof15} and references therein) and $\ubar=0$, then we have
$T_\alpha^\delta(u)=T_\alpha^\delta(-u)$,
which illustrates the non-uniqueness phenomenon.
On the other hand, it should be mentioned that the properties of \tikreg in
Hilbert spaces are well investigated when the penalty functional in the Tikhonov
functional is replaced by $ \myu \mapsto \normix{\myu-\ubar}^2 $, cf., e.g.,
\cite[Chapter 10]{EHN96} or \cite[Section 3.1]{Scherzetal09} and the
references therein, respectively.
}\end{remark}

One of the two main goals of this study is to discuss convergence results for the \tikreg with oversmoothing penalty, i.e.~$ \ust \not \in \ix_1 $ (note, however, that this is not explicitly required anywhere), and the regularization error $ \upardel - \ust $ is still measured in the norm of $ \ix $.
This continues former studies like \cite{GHH19} under the assumption
$ \ust \in \ix_p $ for some $ 0 < p < 1 $. In contrast to those papers, the focus
of the present work is, although also not explicitly required anywhere,
on the case $ \ust \not \in \ix_p $ for each $ 0 < p < 1 $, consequently on the situation characterized by $p=0$. On the other hand, we also mention convergence assertions for $\ust \in X_p$ with $p \ge 1$ under the inequality chain \eqref{eq:normequiv}.

\section{Auxiliary elements and preparatory results}
\label{auxels}
\subsection{Auxiliary elements}
In this section, we consider {\sl auxiliary elements}, which are needed to verify our convergence results. As a preparation, we introduce the bounded, injective, selfadjoint, positive semidefinite linear operator
\begin{align} \label{eq:G}
\Hop \defeq B^{-(\fab)}: \ix \to \ix,
\end{align}
where the operator $ B $ obeying the condition \eqref {eq:b-inverse-stable} is defined in Section~\ref{intro}, and $ a > 0 $ is introduced by item (f) of Assumption~\ref{th:main_assump}. Note that the range $\mathcal{R}(G)$
of $G$ is not closed, and hence zero is an accumulation point of the spectrum $\sigma(G) \subset [0,\|G\|]$ of $G$. In this context, we also mention that $u \in X_p\;(p>0)$ is equivalent to $u \in \mathcal{R}(G^{\frac{p}{2a+2}})$, which means that $u$ obeys a power-type source condition $u=G^{\frac{p}{2a+2}}w$ with some source element $w\in X$. In the case $p=0$, i.e.~if $u \in X$, but $u \notin X_p$ for all $p>0$, then it was shown in \cite{MatHof08} and \cite{HMW09} that there exist an index function\footnote{According to \cite{MatPer03}, we call a function $\varphi\colon (0,\infty) \to (0,\infty)$ index function, if it is continuous, \nondecreasing and satisfies the limit condition $\lim_{t \to 0}\varphi(t) = 0$.}$\varphi$ (for example of logarithmic type, cf.~\cite{Hohage00}) and a source element $w\in X$ such that a (low order) source condition $u=\varphi(G)w$ is satisfied.

The auxiliary elements based on the operator $G$ from \eqref{eq:G} are defined as follows:
\begin{align}
\ualpaux \defeq \ubar + \Hop (\Hop + \para I)^{-1} (\ust-\ubar)
=
\ust - \para(\Hop + \para I)^{-1} (\ust-\ubar)
\for \para > 0,
\label{eq:uaux-def}
\end{align}
where the solution $ \ust$ of the operator equation \eqref{eq:opeq} and the corresponding initial guess $\ubar $ are as introduced above.
The basic properties of the auxiliary elements are summarized in Lemma~\ref{th:auxel}.

\bh{We should mention that the auxiliary elements $\ualpaux$ are the uniquely determined minimizers of the {\sl artificial Tikhonov functional}
$$T_{a,\alpha}(u):=\|u-\ust\|^2_{-a}+ \alpha\,\|u-\ubar\|_1^2$$
over all $u \in X$. The mapping $\ust \mapsto \ualpaux$ is a variant of a {\sl proximal operator} and possesses an explicit character. This allows for error estimates and convergence assertions
for our Hilbert scale model in the case of oversmoothing penalties. If we leave the Hilbert scales, it becomes much more difficult to handle oversmoothing penalties. This is exemplified by the work \cite{GH19},
where the $\ell^1$-regularization is studied when the solution $\ust$ is only in $\ell^2$. There, the auxiliary elements are constructed by projection mappings instead of proximal mappings, and the
occurring conditions for convergence are difficult to interpret.}

In order to specify the limit behaviour of different positive functions occurring in error estimates, we use in the sequel a collection of non-negative functions named  $f_i(\para)\;(i=1,2,\ldots)$ and defined for $ \alpha >0 $ with the property
\begin{align}
\label{eq:fi}
\lim_{\para \to 0} f_i(\para)=0,
\end{align}
to be supposed for all indices $i$.
Consequently, we have for all $i$ that $f_i(\para)=o(1)$ as $ \para \to 0$. Note that pairwise
products $f_i(\para)f_j(\para)$ and linear combinations $K_if_i(\para)+K_jf_j(\para)$ with non-negative constants $K_i,K_j$ can again be written as such a function $f_k(\para)=o(1)$ as $ \para \to 0.$
\begin{lemma}
\label{th:auxel}
There are functions $f_i(\alpha)\;(i=1,2,3)$ for $\alpha>0$
satisfying \refeq{fi}
such that the auxiliary elements from \refeq{uaux-def} have the following properties:
\begin{myenumerate}
\item \label{it:auxel-a}
 $ \normix{\ualpaux - \ust}=f_1(\para) = o(1)\; $ as $\; \para \to 0 $,
\item \label{it:auxel-b}
$ \norma{\ualpaux - \ust} = f_2(\para)\,\para^{\lfrac{a}{(\fab)}} = \landauno{\para^{\lfrac{a}{(\fab)}}} \;$ as $\; \para \to 0 $,
\item \label{it:auxel-c}
$ \normone{\ualpaux-\ubar} = f_3(\para)\,\para^{-\lfrac{1}{(\fab)}}  =\landauno{\para^{-\lfrac{1}{(\fab)}}}\; $ as $ \;\para \to 0 $.
\end{myenumerate}
\end{lemma}
\proof
We show first that
\begin{align}
\label{eq:auxel-0}
\normix{\Hop^\mytheta (\Hop + \para I)^{-1} \myu }
= \landauno{\para^{\mytheta-1}} \as \para \to 0
\end{align}
holds for all $0 \le \mytheta < 1$ and $\myu \in \ix$. It is well known that
\begin{align*}
\norm{\Hop(\Hop + \para I)^{-1} } \le 1,
\quad
\norm{(\Hop + \para I)^{-1} } \le
\para^{-1}
\for \para > 0.
\end{align*}
Then the interpolation inequality implies the estimate
\begin{align}
\label{eq:auxel-d}
 \norm{\Hop^\mytheta (\Hop + \para I)^{-1} }
\le \para^{\mytheta-1}
\for \para > 0 \quad \mbox{and}\quad 0 \le \mytheta \le 1.
\end{align}
Note that the operator $ \Hop $ is selfadjoint and positive semidefinite, and thus
the fractional powers $ \Hop^\mytheta $ are well-defined.

In addition, for fixed $ 0 \le \theta < 1 $ and any $ \myu \in \R(\Hop^q) $ with $ q > 0 $ chosen so small such that
$ \theta + q \le 1 $, we have,
from \refeq{auxel-d} with $ \theta $ replaced by $ \theta + q $,
\begin{align}
\label{eq:auxel-e}
\para^{1-\mytheta} \normix{\Hop^\mytheta (\Hop + \para I)^{-1} \myu }
=
\para^{1-\mytheta} \normix{\Hop^\mytheta (\Hop + \para I)^{-1} \Hop^q v }
\le
\para^{q} \normix{ v } \to 0 \assh \para \to 0,
\end{align}
where $ \myu = \Hop^q v $.
\rp{
The asymptotics \refeq{auxel-0} now follows
from \refeq{auxel-d} and \refeq{auxel-e} and from
an application of the Banach--Steinhaus theorem
to the operators
$ \para^{1-\mytheta}\Hop^\mytheta (\Hop + \para I)^{-1} $
for $ \para \to 0 $.
Here, $0 \le \theta < 1$ is fixed, and
we have used the fact that for arbitrary $ q > 0 $, the range of the operator $ \Hop^q $ is dense in $ \ix $, i.e.~$ \overline{\R(\Hop^q)} = \ix $.
For the Banach--Steinhaus theorem,
cf., e.g., \cite[Problem~10.1]{Kress[89]} or \cite[Theorem 1.1.4]{Saranen_Vainikko[01]}.
}

For the functions
\begin{align} \label{eq:f1}
f_1(\para) &= \|\para(\Hop + \para I)^{-1} (\ust-\ubar)\|,
\\
\label{eq:f2}
f_2(\para) &=
\para^{-\lfrac{a}{(\fab)}} \| \Hop^{\lfrac{a}{(\fab)}} [\para (\Hop + \para I)^{-1} (\ust-\ubar) ] \|, \\
\label{eq:f3}
f_3(\para) &=
\para^{\lfrac{-(2a+1)}{(\fab)}}\|\Hop^{\lfrac{(2a+1)}{(\fab)}}
[\para (\Hop + \para I)^{-1} (\ust - \ubar) ]\|,
\end{align}
the statements of the lemma are now easily obtained from
\refeq{auxel-0} and the following three representations,
\begin{align*}
\ualpaux - \ust & = - \para(\Hop + \para I)^{-1} (\ust-\ubar),
\\
B^{-a} ( \ualpaux - \ust)
&=- \Hop^{\lfrac{a}{(\fab)}} [\para (\Hop + \para I)^{-1} (\ust-\ubar) ],
\\
B (\ualpaux-\ubar) &= \Hop^{\lfrac{(2a+1)}{(\fab)}} [(\Hop + \para I)^{-1} (\ust - \ubar)].
\proofend
\end{align*}
\subsection{Some estimates for oversmoothing Tikhonov regularization}
$\;$

\begin{lemma}
\label{th:upardel_lemma}
\mainassump
Then there is a function $f_4(\para)$ for $\para>0$ satisfying \refeq{fi}
such that for all \rpc{$ \para > 0 $ and} $\delta>0$,
\rpc{we have}
\begin{align*}
\max\{\normyps{\F{\upardel} - \fdel }, \, \sqrt{\para} \normone{\upardel-\ubar}\}
\le f_4(\para)\,\para^{\lfrac{a}{(\fab)}}+ \delta.
\end{align*}
\end{lemma}
\proof For $ \para > 0 $ small enough, {\rp{say $ 0 < \para \le \para_0 $,} we have
$ \ualpaux \in \Dset $, because item (a) of Lemma~\ref{th:auxel} holds and $\ust$ is an interior point of $\DF$. Thus
\begin{align*}
& \kla{\normypsqua{\F{\upardel} - \fdel } + \para \normonequa{\upardel-\ubar}}^{1/2}
\le \kla{\normqua{\F{\ualpaux} - \fdel } + \para \normonequa{\ualpaux-\ubar}}^{1/2}
\\
& \quad \le \normyps{\F{\ualpaux} - \fdel } + \sqrt{\para} \normone{\ualpaux-\ubar}
\le \normyps{\F{\ualpaux} - \fst } + \sqrt{\para} \normone{\ualpaux-\ubar} + \delta.
\end{align*}
The first term on the \rhs of the latter estimate can be written as
\begin{align*}
\normyps{\F{\ualpaux} - \fst }
\le \cb \norma{\ualpaux - \ust}
\le \cb \,f_2(\para)\,\para^{\lfrac{a}{(\fab)}}.
\end{align*}
This is
a consequence of item \ref{it:auxel-b} of Lemma \ref{th:auxel}.
The second term on the \rhs of the latter estimate attains the form
\begin{align*}
\sqrt{\para} \normone{\ualpaux-\ubar}\le f_3(\para)\,\para^{\lfrac{a}{(\fab)}},
\end{align*}
based on item \ref{it:auxel-c} of Lemma \ref{th:auxel}. This yields the function
\rp{
$$f_4(\para):=  \cb \,f_2(\para)+f_3(\para) \quad \mbox{for} \;\para \le \para_0. $$
Note that $f_4(\para) \to 0 $ as $ \para \to 0. $
For $ \para > \para_0 $, the estimate
\begin{align*}
\kla{\normypsqua{\F{\upardel} - \fdel } + \para \normonequa{\upardel-\ubar}}^{1/2}
\le \normyps{\F{\ubar} - \fdel } \le
\normyps{\F{\ubar} - \fst } + \delta,
\end{align*}
holds, so we may define
$$f_4(\para):=
\frac{\normyps{\F{\ubar} - \fst }}{\para^{\lfrac{a}{(\fab)}}}
\quad \mbox{for} \;\para > \para_0. $$
This completes the proof of the lemma.
}
\proofend
\begin{corollary}
\label{th:upardel_cor}
\mainassump
Then there \rp{are} a function $f_5(\para)$ for $\para>0$ satisfying \refeq{fi} and a constant $K_1>0$ such that for all \rpc{$ \para > 0 $ and} $\delta>0$, \rpc{we have}
$$ \norma{\upardel -\ust} \le f_5(\para)\,\para^{\lfrac{a}{(\fab)}} + K_1 \,\delta.$$
\end{corollary}
\proof
It follows from the left-hand estimate in \refeq{normequiv} and Lemma \ref{th:upardel_lemma} that
\begin{align*}
& \ca \norma{\upardel -\ust} \le \normyps{\F{\upardel} - \fst } \le
\normyps{\F{\upardel} - \fdel } + \delta
\le f_4(\para)\,\para^{\lfrac{a}{(\fab)}} + 2\delta.
\end{align*}
\rpc{The assertion of the corollary now follows by setting $f_5(\para):=\tfrac{f_4(\para)}{\ca}$ and $K_1:=\tfrac{2}{\ca}.$}
\proofend

\bn
The error $ \norm{\upardel -\ust} $ is now estimated by the following series of error estimates.
Using the triangle inequality and Lemma~\ref{th:auxel}, we obtain
\begin{align}
\normix{\upardel -\ust}\le \normix{\upardel -\ualpaux}+\normix{\ualpaux-\ust}=\normix{\upardel -\ualpaux}+f_1(\para),
\label{eq:upardel_prop-a}
\end{align}
and below we consider the term $ \normix{\upardel -\ualpaux} $ in more detail.
From the interpolation inequality for bounded linear, self-adjoint and positive semidefinite operators on Hilbert spaces, cf.~\cite[(2.49)]{EHN96}, it follows
\begin{align}
\label{eq:upardel_prop-b}
\normix{\upardel -\ualpaux} & \le
\norma{\upardel -\ualpaux}^{\lfrac{1}{(a+1)}}
\normone{\upardel -\ualpaux}^{\lfrac{a}{(a+1)}}.
\end{align}
Both terms on the \rhs of the  estimate \eqref{eq:upardel_prop-b} can be estimated by using Corollary~\ref{th:upardel_cor} and Lemma~\ref{th:auxel} in the following manner:
Precisely, we find with $f_6(\para):=f_2(\para)+f_5(\para)$ and $f_7(\para):=f_3(\para)+f_4(\para)$ the estimates
\begin{align*}
\norma{\upardel -\ualpaux}
& \le
\norma{\upardel -\ust} + \norma{\ualpaux - \ust} \le f_6(\alpha)\,\para^{\lfrac{a}{(\fab)}}+K_1\delta,\\
\normone{\upardel -\ualpaux}
& \le
\normone{\upardel-\ubar} + \normone{\ualpaux-\ubar}
\le  \para^{-\lfrac{1}{2}}\left(f_7(\para)\, \para^{\lfrac{a}{(\fab)}}+\delta\right).
\end{align*}
Thus we can continue estimating \refeq{upardel_prop-b}.
Introducing $ f_8(\para):= \max\{ f_6(\para), f_7(\para) \} $
and $ K_2 := \max\{K_1, 1\} $, we obtain
\begin{align*}
\normix{\upardel -\ualpaux}
& \le \left(f_6(\alpha)\,\para^{\lfrac{a}{(\fab)}}+K_1\delta \right)^{1/(a+1)}\,\left(\para^{-\lfrac{1}{2}}\left(f_7(\para)\, \para^{\lfrac{a}{(\fab)}}+\delta\right) \right)^{a/(a+1)}
\\
& \le \left(f_8(\alpha)\,\para^{\lfrac{a}{(\fab)}}+K_2\delta \right)^{1/(a+1)}\,\left(\para^{-\lfrac{1}{2}}\left(f_8(\para)\, \para^{\lfrac{a}{(\fab)}}+ K_2\delta\right) \right)^{a/(a+1)}
\\
& =
\para^{\lfrac{-a}{(\fab)}}
 \left(f_8(\alpha)\,\para^{\lfrac{a}{(\fab)}}+K_2\delta \right)
 =
f_8(\alpha) + K_2
\frac{\delta}{\para^{\lfrac{a}{(\fab)}}}.
\end{align*}
From the latter estimate and \refeq{upardel_prop-a}, the following proposition now immediately follows by considering $ f_9(\para) := f_1(\para) + f_8(\para) $ there.
\begin{proposition}
\label{th:upardel_prop}
\mainassump
Then there \rp{are} a function $f_9(\para)$ for $\para>0$ satisfying \refeq{fi} and a constant $K_2>0$ such that for all \rpc{$ \para > 0 $ and} $\delta>0$,
\rpc{we have}
\begin{align}
\label{eq:fin}
\norm{\upardel -\ust} \le f_9(\para) + K_2 \frac{\delta}{\para^{\lfrac{a}{(\fab)}}}.
\end{align}
\end{proposition}
The inequality \eqref{eq:fin}, which is valid for arbitrary noise levels $\delta>0$ and regularization parameters $\para>0$, allows us to formulate in the subsequent section sufficient conditions for the convergence of the error norm $\normix{\upardel -\ust}$ of the regularized solutions.
%
\section{Convergence results}
\label{parameter_choices}
\subsection{Main theorem}
The following main theorem is an immediate consequence of the error estimates outlined in the preceding section. The formulated convergence result follows immediately from the inequality~\eqref{eq:fin}.
\begin{theorem}
\label{th:main}
\mainassump
Then for any a priori parameter choice $\para_*=\para(\delta)$ and any a posteriori parameter choice  $\para_*=\para(\delta,y^\delta)$, the regularized solutions $u_{\para_*}^\delta$ converge in the norm of the Hilbert space $X$ to the solution $\ust$ of the operator equation \eqref{eq:opeq} for $\delta \to 0$, i.e.~$\lim \nolimits_{\delta \to 0}\normix{u_{\para_*}^\delta -\ust}=0$, whenever
\begin{equation} \label{eq:mainthm}
\para_* \to 0 \qquad \mbox{and} \qquad \mfrac{\delta}{\para_*^{\lfrac{a}{(\fab)}}} \to 0 \qquad \mbox{as} \qquad \delta \to 0.
\end{equation}
\end{theorem}

\begin{remark} \label{rem:rem3}
If\footnote{\rp{With a slight abuse of notation, for two nonnegative functions we write $ f(\delta) \sim g(\delta) $, if there two constants $ 0 < c_1 \le c_2 $ such that $ c_1 f(\delta) \le g(\delta) \le c_2 f(\delta) $ for each $ \delta > 0 $ small enough.}}%
\begin{equation} \label{eq:borderline}
\para_*^{\lfrac{a}{(\fab)}} \sim \delta
\qquad \mbox{as} \qquad \delta \to 0,
\end{equation}
then convergence cannot be derived in this way, because in that borderline case, the second term on the right-hand side of inequality \eqref{eq:fin} does not tend to zero.
\end{remark}
\begin{remark} \label{rem:rem4}
\bh{
The convergence result of Theorem~\ref{th:main} applies both to (a) the classical case $\ust \in \DB$ as well as to (b) oversmoothing penalties $\ust \notin \DB$.
This theorem is \rp{directly} based on formula \eqref{eq:fin}. An inspection of the proof of this formula, given by means of Lemma~\ref{th:upardel_lemma} through Proposition~\ref{th:upardel_prop}, shows that both inequalities in \eqref{eq:normequiv}
are needed for the convergence result of Theorem~\ref{th:main}.
Nevertheless, for the oversmoothing case (b), this is real progress, since the convergence result of Theorem~\ref{th:main} does not use any form of additional smoothness on $\ust$. Such {\sl pure convergence} assertions for the oversmoothing case, without any smoothness assumptions like the condition $\ust \in X_p$ for some $p \in (0,1)$ that was used in~\cite{Hofmann_Mathe[18]} for obtaining {\sl convergence rates} results,
are missing in the literature by now.}

As is well known, in case (a) with $\ust \in \DB $, the parameter choice condition
$$\para_* \to 0 \qquad \mbox{and} \qquad \mfrac{\delta^2}{\para_*} \to 0 \qquad \mbox{as} \qquad \delta \to 0,$$
which is stronger than  \eqref{eq:mainthm}, is always sufficient for convergence of regularized solutions, and inequalities occurring in \eqref{eq:normequiv} represent only tools for obtaining convergence rates. On the other hand, in the limit situation $\para_* \sim \delta^2$ of choosing the regularization parameter, one needs the left-hand inequality in \eqref{eq:normequiv} for obtaining convergence rates, and this inequality occurs here as a conditional stability estimate (cf.~\cite[Prop.~3]{GHH19}, \cite[Theorem~1.1]{Egger_Hofmann[18]} and references therein). Convergence is then a consequence of derived convergence rates.
\end{remark}
\subsection{A priori parameter choice of power type}
In this subsection, we consider in light of Theorem~\ref{th:main} the a priori parameter choice
\begin{equation} \label{eq:monomial}
\para_*=\para(\delta) \sim  \delta^\kappa
\end{equation}
for exponents $\kappa>0$. Then condition \eqref{eq:mainthm} is satisfied if and only if $0<\kappa < 2+\frac{2}{a}$, and the borderline condition \eqref{eq:borderline} holds if and only if $\kappa=2+\frac{2}{a}$. This gives the following proposition.
\begin{proposition}\label{pro:monomial}
For the a priori choice \eqref{eq:monomial} of the regularization parameter $\para>0$, condition \refeq{mainthm} in Theorem~\ref{th:main}
holds if and only if $0<\kappa < 2+\frac{2}{a}$. For all $a>0$, the choice $\para_* \sim \delta^2$ yields convergence.
\end{proposition}

We can distinguish the $\kappa$-intervals \,(A): $0< \kappa < 2$,\, (B): $\kappa=2$, and  (C): $2<\kappa < 2+\frac{2}{a}$ for \eqref{eq:monomial}. Then we have $\frac{\delta^2}{\alpha_*} \to 0$ as $\delta \to 0$ in situation (A) and $\frac{\delta^2}{\alpha_*} \sim 1$
in situation (B). Note that both situation also occur and yield convergent regularized solutions in the oversmoothing case $\ust \notin \DB$. This is a bit surprising, because the behaviour
\begin{equation} \label{eq:infty}
\frac{\delta^2}{\alpha_*} \to \infty \quad \mbox{as} \quad \delta \to 0
\end{equation}
occurring in situation (C) was supposed in the literature to be typical for the case of oversmoothing penalties. Namely as is seen in \cite{Hofmann_Mathe[19]}, convergence rate results of the form
$$\normix{u_{\para_*}^\delta -\ust}= \mathcal{O}\left(\delta^{\frac{p}{a+p}} \right) \quad \mbox{as} \quad \delta \to 0 $$
are obtained under the both-sided structural condition \eqref{eq:normequiv} and in particular under the smoothness assumption
$\ust \in X_p$ for $0<p<1$, whenever the a priori parameter choice of type \eqref{eq:monomial} with prescribed exponent $\kappa= \frac{2(a+1)}{a+p}=2+\frac{2(1-p)}{a+p}$ applies. Evidently, this prescribed $\kappa$ satisfies
the conditions \eqref{eq:infty} and $2<\kappa< 2 + \frac{2}{a}$ for all $0<p<1$. It is important to note that $p=0$ coincides with the borderline case $\kappa=2+\frac{2}{a}$ which, however, is not sufficient for convergent regularized solutions.
\subsection{A discrepancy principle}
For the specification of an appropriate discrepancy principle,
the behaviour of the misfit functional $ \para \mapsto \Jdel{\upardel} $ needs to be described, for $ \delta > 0 $ fixed.
The basic properties are summarized in the following proposition.
\begin{proposition}
\label{th:misfit-behavior}
\mainassump
Then for $ \delta > 0 $ fixed, the function $ \para \mapsto \Jdel{\upardel} $ is \monotonically \nondecreasing, with
\begin{align}
\lim_{\para \to 0}  \Jdel{\upardel} \le \delta,
\quad \lim_{\para \to \infty}  \Jdel{\upardel} =  \Jdel{\ubar}.
\label{eq:misfit-behavior}
\end{align}
We have $ \lim_{\para \to \infty} \normix{\upardel - \ubar} = 0 $.
\end{proposition}
\proof
We start with the verification of the first statement of the proposition. As a preparation, we show that
the function $ \para \mapsto \normone{\upardel-\ubar} $ is \monotonically non-increasing. Indeed, for $ 0 < \para \le \parac $ fixed,
we have
\begin{align*}
\pad[\parac]{\upardelb} & \le \pad[\parac]{\upardel}
= \pad{\upardel} + (\parac-\para) \normonequa{\upardel-\ubar}
\\
& \le \pad{\upardelb} + (\parac-\para) \normonequa{\upardel-\ubar}
\\
& = \pad[\parac]{\upardelb} + (\parac-\para) (\normonequa{\upardel-\ubar} - \normonequa{\upardelb-\ubar}),
\end{align*}
and thus
$ \normone{\upardelb-\ubar} \le \normone{\upardel-\ubar} $.
The first statement of the proposition is now easily obtained:
for $ 0 < \para \le \parac $ we have
\begin{align*}
\Jdelqua{\upardel} + \para \normonequa{\upardel-\ubar}
& = \pad{\upardel} \le \pad{\upardelb} =
 \Jdelqua{\upardelb} + \para \normonequa{\upardelb-\ubar}
\\ &
\le \Jdelqua{\upardelb} + \para \normonequa{\upardel-\ubar},
\end{align*}
and thus
$ \Jdel{\upardel} \le \Jdel{\upardelb} $.

Next we consider the latter statement of the proposition. There holds
\begin{align}
\label{eq:misfit-behavior-a}
\Jdelqua{\upardel} + \para \normonequa{\upardel - \ubar} = \pad {\upardel} \le \pad {\ubar} =  \Jdelqua{\ubar},
\end{align}
and thus in particular
$ \normone{\upardel - \ubar} = \Landauno{\para^{-1/2}} $ as $ \para \to \infty $.
The estimate \refeq{b-inverse-stable}
implies $ \normix{\upardel - \ubar} = \Landauno{\para^{-1/2}} $ as $ \para \to \infty $, which implies the latter statement of the proposition.

The first statement in \refeq{misfit-behavior} follows directly from
Lemma \ref{th:upardel_lemma},
and we finally consider the second statement in \refeq{misfit-behavior}.
From \refeq{misfit-behavior-a} we already know that
$ \lim_{\para \to \infty}  \Jdel{\upardel}$ $\le  \Jdel{\ubar} $.
Conversely, sequential weak continuity of the operator $ \Fpur $ implies weak convergence
$ \F{\upardel} \rightharpoonup \F{\ubar} $ as $ \para \to \infty $, and thus
$ \lim_{\para \to \infty}  \Jdel{\upardel} \ge  \Jdel{\ubar} $. This completes the proof of the proposition.
\proofend
\begin{remark}
Notice that in the proof of Proposition \ref{th:misfit-behavior},
no use of the fundamental estimates \refeq{normequiv} on the smoothing property of $ F $
is made in fact.
Notice also that the statement of Proposition \ref{th:misfit-behavior} is quite similar to related results for \tikreg with non-oversmoothing penalty,
cf.~\cite{Anzengruber_Hofmann_Mathe[14]},
\cite[Section 2.6]{Tikhonov_Leonov_Yagola[98]} and \cite[Section 6.7]{Vainikko[91.1]}.
\end{remark}

It follows from Proposition \ref{th:misfit-behavior} that the following version of the discrepancy principle (cf.~\cite{Vainikko[88],Vainikko[91.1]}) is implementable.
It determines, for each noise level $ \delta > 0 $, an approximation $ \udel \in \Dset $.
Possibly discontinuities of the misfit functional $ \para \mapsto \Jdel{\upardel} $
are taken into account.
\begin{algorithm}[Discrepancy principle]
\label{th:discrepancy-def}
Let $ \myb > 1 $ and $ \myc > 1 $ be finite constants.
\begin{myenumerate}
\item
If $ \Jdel{\ubar} \le \myb \delta $ holds,
then choose $ \pardel = \infty $, i.e.~$ \upardel[\infty] \defeq \ubar \in \Dset $.

\item
Otherwise choose a finite parameter $ \para =: \pardel > 0 $ such that
\begin{align}
\Jdel{\updd} \le \myb \delta \le \Jdel{\upddb}
\quad \text{for some } \pardel \le \pardelb \le \myc \pardel.
\label{eq:parameter-choice}
\end{align}
\end{myenumerate}
\end{algorithm}
Algorithm \ref{th:discrepancy-def} can be realized by the following strategy.
\begin{remark}[Sequential discrepancy principle]
\label{th:sdp}
Practically, a parameter $ \pardel $ satisfying condition \refeq{parameter-choice}
can be determined, e.g., by choosing a constant $ \theta > 1 $ and an initial guess
$ \parainit > 0 $, and proceeding then as follows:
\begin{mylist}
\item
If $ \Jdel{\uadno[0]} \ge \myb \delta $ holds,
then, with the notation $ \para^{(k)} = \theta^{-k} \parainit $,
proceed for $ k = 1, 2, \ldots $ until
$ \Jdel{\uadno} \le \myb \delta \le \Jdel{\uadno[k-1]} $
is satisfied for the first time;
define $ \pardel = \alpha^{(k)} $ then.
\item
If $ \Jdel{\uadno[0]} \le \myb \delta $ holds,
then, with the notation $ \alpha^{(k)} = \theta^{k} \parainit $,
proceed for $ k = 1, 2, \ldots $ until
$ \Jdel{\uadno[k-1]} \le \myb \delta \le \Jdel{\uadno} $
is satisfied for the first time;
define $ \pardel = \para^{(k-1)} $ then.
\end{mylist}
\end{remark}

The regularizing properties of Algorithm \ref{th:discrepancy-def}
are stated in the following theorem.
\begin{theorem}
\label{th:aposteriori-convergence}
\mainassump
For the a posteriori parameter choice introduced in Algorithm \ref{th:discrepancy-def},
we have
\begin{align}
\normix{\udel -\ust} \to 0, \quad
\tfrac{\delta}{\pardel^{\lfrac{a}{(\fab)}}} \to 0
\as \delta \to 0.
\label{eq:aposteriori-convergence}
\end{align}
\end{theorem}
\proof For an arbitrary countable noise level set $ \myH \subset \reza_+ $ having the origin as only accumulation point, we consider the following three cases:
(a) $ \pardel = \infty $ for each $ \delta \in \myH $,
(b) $ \pardel \to 0 $ as $ \myH \ni \delta \to 0 $, and
(c) $ \pardel < \infty $ for each $ \delta \in \myH $,
$ \liminf_{\myH \ni \delta \to 0} \pardel > 0 $. Below we show that in each of those three cases, \refeq{aposteriori-convergence} holds.
The main statement of the theorem then follows by arguing for subsequences.
Note that in cases (a) and (c), the second statement in \refeq{aposteriori-convergence} trivially holds.
\begin{myenumerate}
\item
The case $ \pardel =  \infty $ for $ \delta \in \myH $ means
$ \Jdel{\ubar} \le \myb \delta $ for $ \delta \in \myH $, and thus
$ \F{\ubar} = \fst $, and then $ \upardel[\infty] = \ubar = \ust $ for
$ \delta \in \myH $.

\item
Suppose that $ \pardel \to 0 $ as $ \myH \ni \delta \to 0 $. From Lemma \ref{th:upardel_lemma}, we obtain
\begin{align*}
\myb \delta \le \normyps{\F{\upddb} - \fdel }
\le \landauno{\pardelb^{\lfrac{a}{(\fab)}}} + \delta
= \landauno{\pardel^{\lfrac{a}{(\fab)}}} + \delta,
\end{align*}
and thus
$ \lfrac{\delta}{\pardel^{\lfrac{a}{(\fab)}}} \to 0 $
as $ \myH \ni \delta \to 0 $.
Proposition \ref{th:upardel_prop}
then yields $ \normix{\updd -\ust} \to 0 $ as $ \myH \ni \delta \to 0 $.

\item
Next suppose that both
$ \pardel < \infty $ for $ \delta \in \myH $, and
\begin{align}
\label{eq:pardel_not_to_zero}
\liminf_{\myH \ni \delta \to 0} \pardel > 0.
\end{align}
\begin{myenumerate_roman}
\item
We first observe that
\begin{align}
\label{eq:pardel_not_to_zero_a}
\ca\norma{\updd - \ust} \le \normyps{\F{\updd} - \fst} \le
\normyps{\F{\updd} - \fdelta} + \delta \le (\myb +1 ) \delta,
\end{align}
so
$ \norma{\updd - \ust} = \Landauno{\delta} $ as $ \myH \ni \delta \to 0 $.
Note that the asymptotics \refeq{pardel_not_to_zero} is not needed for this result.

\item
There holds
\begin{align}
\label{eq:u-DB-a}
\normone{\updd} = \Landauno{1} \as \myH \ni \delta \to 0.
\end{align}
This easily follows from \refeq{misfit-behavior-a} and
\refeq{pardel_not_to_zero}, in combination with the estimate
$ \sqrt{\pardel} \normone{\updd - \ubar} \le \Jdel{\ubar}
\le \J{\ubar} + \delta $.

\item
We next show
\begin{align}
\label{eq:u-DB-b}
\ust \in \DB.
\end{align}
\rp{For this purpose, we observe that estimate \refeq{u-DB-a}
implies
weak convergence in $ \ix_1 $ for some subsequence $ \myHp \subset \myH $, \ie
for some element $ \myv \in \DB = \ix_1 $, we have
$ \updd \rightharpoonup \myv $ in $ \ix_1 $ as $ \myHp \ni \delta \to 0 $.
From the (weak) continuity of the embedding operator $ \ix_1 \hookrightarrow \ix $, we then obtain $ \updd \rightharpoonup \myv $ in $ \ix $ as $ \myHp \ni \delta \to 0 $, and thus
$ \myv \in \Dset $, due to the weak closedness of $ \DF $.
Since the operator $ F $ is sequentially weakly continuous, we have
$ \F{\updd} \rightharpoonup \F{\myv} $ as $ \myHp \ni \delta \to 0 $.
Algorithm \ref{th:discrepancy-def} implies $ \normyps{\F{\updd} - \F{\ust}} \to 0 $ as $ \delta \to 0 $, so we finally obtain
$ \F{\myv} = \F{\ust} $. The lower bound in \refeq{normequiv}
then gives $ \myv = \ust $, which finally implies
\refeq{u-DB-b}.}
\item
From \refeq{pardel_not_to_zero_a}, \refeq{u-DB-b} and the interpolation inequality, we now obtain
\begin{align*}
\normix{\updd -\ust} & \le
\norma{\updd -\ust}^{\lfrac{1}{(a+1)}} \cdot
\normone{\updd -\ust}^{\lfrac{a}{(a+1)}}
=
\Landauno{\delta^{\lfrac{1}{(a+1)}}} \as \myH \ni \delta \to 0.
\end{align*}
This completes the proof of the theorem.
\proofend
\end{myenumerate_roman}
\end{myenumerate}
Note that in the oversmoothing situation $ \ust \not \in \ix_1 $, case (c) in the proof of Theorem~\ref{th:aposteriori-convergence} does not emerge, cf.~\refeq{u-DB-b}. This fact is, in a similar setting, already observed in \cite[Lemma~1]{Hofmann_Mathe[18]}.
\begin{remark}
Notice that the situation (b) in the proof of Theorem~\ref{th:aposteriori-convergence}
 is the regular case in applications.  The case (c) is an exceptional case which,
in the non-oversmoothing case, can be excluded,
if the exact penalization veto is satisfied. This veto had been introduced in the paper \cite{Anzengruber_Hofmann_Mathe[14]}; see also \cite{Anzengruber_ramlau[10]}.
\end{remark}
\section{Low order convergence rates} \label{sec:lowrates}
Our convergence assertion established in the main theorem formulated in Subsection~\ref{parameter_choices}, is due to the error estimate \eqref{eq:fin} derived in Section~\ref{auxels}. The presented sufficient conditions
for convergence are based on the Banach--Steinhaus theorem and do not need any form of solution smoothness. In other words, the case $p=0$ is included, where $\ust$ does not satisfy a power-type source condition. However,
as already mentioned above, there exists at least a source condition of lower order for solution element $\ust \in X$. Precisely, there is always an index function $\varphi$ and a source element $w \in X$ such that
\begin{equation} \label{eq:losc}
\ust-\ubar = \varphi(G)\, w.
\end{equation}
Based on formula \eqref{eq:fin} and taking into account the representations \eqref{eq:f1}, \eqref{eq:f2} and \eqref{eq:f3}, we can derive for such source condition low order convergence rates in the case of oversmoothing penalties as a byproduct of the studies presented in Section~\ref{auxels}. We will outline this in the following.
\begin{lemma} \label{lem:quali}
\bh{If, for an index function $\varphi$, the quotient function $\varphi(t)/t$ is \nonincreasing for $0<t \le \overline t$ with some constant $\overline t \in (0,\|G\|]$,
then there exist positive constants $C$ and $\overline \alpha$ such that
\begin{equation} \label{eq:qualia}
\sup \limits_{0<\lambda \le \|G\|} \frac{\alpha \varphi(\lambda)}{\lambda+\alpha} \le C\,\varphi(\alpha) \qquad (0<\alpha \le \overline \alpha).
\end{equation}
}\end{lemma}
\bh{The assertion of the lemma follows directly from \cite[Prop.~3.3]{BHTY06}.\footnote{\bh{For the understanding of the formula
(\ref{eq:qualia}), the concept of  {\sl qualification} for a regularization method introduced in \cite{MatPer03} is helpful. Precisely, all index functions $\varphi(t)$ which are covered by the function $t$ are qualifications for the classical Tikhonov regularization applied to the operator $G$, which implies that
an inequality of type (\ref{eq:qualia}) is valid.}}}
\begin{corollary} \label{cor:quali}
\rp{
Let $\varphi$ be an index function such that for each exponent $\eta>0$ the quotient function $t^\eta/\varphi(t)$ is \strictly increasing for sufficiently small $t>0$.
Then for each $ 0 \le \theta < 1 $} \bh{there exist positive constants $C$ and $\overline \alpha$ such that
\begin{equation} \label{eq:theta}
\sup \limits_{0<\lambda \le \|G\|} \frac{\alpha \lambda^\theta\varphi(\lambda)}{\lambda+\alpha} \le C\,\alpha^\theta\varphi(\alpha) \qquad (0<\alpha \le \overline \alpha).
\end{equation}
}
\end{corollary}
\bh{
\proof
We have that for all $0 \le \theta <1$, the quotient function $\frac{t^\theta \varphi(t)}{t}=\frac{\varphi(t)}{t^{1-\theta}}$ with $1-\theta >0$ is \nonincreasing for sufficiently small $t>0$. Consequently, there are, according to formula (\ref{eq:qualia}) of Lemma~\ref{lem:quali},
positive constants $C$ and $\overline \alpha$ depending on $\theta$ such that (\ref{eq:theta}) is valid.
\proofend
}
\begin{theorem}
\label{th:lowrates}
Let Assumption \ref{th:main_assump}
and the source condition \eqref{eq:losc} be satisfied, where it is supposed that for all $\eta>0$, the index function $\varphi$ has a \strictly increasing quotient function $t^\eta/\varphi(t)$ for sufficiently small $t>0$.
Then we have, for some positive constant $K_0$ and $K_2$ from \eqref{eq:fin}
 and for all $\delta>0$ and sufficiently small $\alpha>0$, the error estimate
\begin{align}
\label{eq:fin2}
\norm{\upardel -\ust} \le K_0\,\varphi(\para) + K_2 \,\frac{\delta}{\para^{\lfrac{a}{(\fab)}}}.
\end{align}
\end{theorem}
\proof
\bh{Based on the source condition (\ref{eq:losc})}, the functions $f_1, f_2$ and $f_3$ from Lemma~\ref{th:auxel} satisfy
\begin{align*}
f_1(\alpha)=\mathcal{O}(\varphi(\alpha)),\quad f_2(\alpha)=\mathcal{O}(\varphi(\alpha)), \quad f_3(\alpha)=\mathcal{O}(\varphi(\alpha)) \quad \mbox{as} \quad \alpha \to 0.
\end{align*}
These properties
are immediate consequences from \eqref{eq:theta}, taking into account the three representations \eqref{eq:f1}, \eqref{eq:f2} and \eqref{eq:f3}.
Since the function $f_9(\para)$ in the error estimate \eqref{eq:fin}
can be estimated from above by linear combinations and maximizing of the functions $f_1, \,f_2$ and $f_3$
\rp{for $ \alpha>0 $ sufficiently small},
there is a positive constant $K_0$ such that
$f_9(\para) \le K_0\,\varphi(\para)$ holds for sufficiently small $\para>0$.
 \proofend

\bigskip \noindent
This provides us directly with the following low order convergence rate result.

\begin{corollary} \label{cor:lowrates}
Set, under the assumptions of \bh{Theorem}~\ref{th:lowrates}, $\psi(\para):=\varphi(\para)\, \para^\frac{a}{2a+2}$ and $\para_*=\para(\delta):=\psi^{-1}(\delta)$. Then we have
$$\|u^\delta_{\para_*} -\ust\|=\mathcal{O}\left( \varphi(\psi^{-1}(\delta) \right) \qquad \textup{as} \quad \delta \to 0.$$
\end{corollary}

\begin{example} \label{ex:logrates}
In this example, we consider source conditions \eqref{eq:losc} of logarithmic type with the function
\begin{align}
\varphi(t)=\varphi_{\log}^\kappa(t):= (-\log(t))^{-\kappa} \qquad (\kappa>0),
\label{eq:log-type-phi}
\end{align}
which is strictly concave for sufficiently small $t>0$, and can be extended to $(0,\infty)$ as an index function. Is is evident for all $\eta,\kappa>0$ that the quotient function
$t^\eta/\varphi_{\log}^\kappa(t)$ is \strictly increasing for sufficiently small $t>0$, and Corollary~\ref{cor:quali} applies. This yields the error estimate \eqref{eq:fin2} written as
$$\norm{\upardel -\ust} \le K_0\,\left(-\log(\para)\right)^{-\kappa} + K_2 \,\frac{\delta}{\para^{\lfrac{a}{(\fab)}}}.$$
For the a priori choice $\para_*=\para(\delta) \sim \delta^2$
 of the regularization parameter, this implies the logarithmic convergence rate
\begin{align}
\|u^\delta_{\para_*} -\ust\|=\mathcal{O}\left( (-\log(\delta))^{-\kappa} \right) =\mathcal{O}\left( \varphi_{\log}^\kappa(\delta) \right)  \qquad \mbox{as} \quad \delta \to 0.
\label{eq:log-rates}
\end{align}
Note that this parameter choice strategy differs from that presented in Corollary \ref{cor:lowrates}.
\end{example}

\section{Numerical illustrations}
The theoretical results are numerically illustrated
for the nonlinear operator $ F : \elltwo \supset \DF \to \elltwo $, given by
the sum $ F = F_1 + F_2 $ of a linear operator $ F_1 $ and a quadratic operator
$ F_2 $ as follows,
\begin{align}
& F_1: \elltwo \supset \DF  \to \elltwo, \quad (u_n) \mapsto 7 (\an u_n),
\label{eq:F1-def} \\
& F_2: \elltwo \supset \DF \to \elltwo, \quad (u_n) \mapsto (\an u_n^2).
\label{eq:F2-def}
\end{align}
Here,
\begin{align*}
\DF =  \inset{ \myu \in \elltwo \mid \normb{\myu} \le 3 },
\end{align*}
and $ \elltwo = \inset{ (u_n) \mid \normbqua{u} = \sum_{n=1}^{\infty} u_n^2 < \infty } $.
The stronger norm $\|\cdot\|_1$ is defined by the generator
\begin{align*}
B: \elltwo \supset \DB \to \elltwo, \quad (u_n) \mapsto (\aninv u_n),
\qquad \DB \defeq \inset{ (u_n) \mid (n u_n) \in \elltwo }.
\end{align*}
In what follows, we consider the equation $ \F u = \fst $ having
\begin{align*}
\ust = (\ust_n), \quad \textup{with }
\ust_1 = 1, \quad
\ust_n =
\frac{1}{\sqrt{n} (\log n)^{\myq}}, \quad n = 2, 3, \ldots,
\end{align*}
as a solution.
Assumption \ref{th:main_assump} is satisfied then; in particular, the two structural inequalities in \refeq{normequiv} are satisfied for $ a = 1 $.
In addition, we have $ \ust \not \in \DB $.

Below, some additional remarks on the numerical tests are given.
\begin{mylist_indent}
\item We consider the \tikreg \refeq{TR} with $ \ubar = 0 $.

\item
For the finite-dimensional approximation needed for the computations, we replace
in
\refeq{F1-def}, \refeq{F2-def}
the space $ \elltwo $ by $ \reza^\N $, with $ \N = 6000 $ at each occurrence.

\item
In the numerical experiments, we consider perturbations of the form
$ f_n^\delta = f_n + \Delta_n $
for $ n = 1,2,\ldots, \N $,
with uniformly distributed random values $ \Delta_n $  satisfying $ \modul{\Delta_n} \le \delta/\sqrt{N} $.
\end{mylist_indent}
For this framework, we consider \tikreg
\refeq{TR} with an a~priori and an a~posteriori parameter choice, respectively.
\subsection{Numerical results for an a priori parameter choice}
\label{example1}
We first consider the a~priori parameter choice $ \pardel = \delta^2 $, for different values of $ \delta $. The numerical results in Table \ref{tab:num1} confirm the logarithmic convergence rate given by \refeq{log-rates}.
Note that a logarithmic type source condition $ \ust = \varphi_{\log}^\kappa(G) w $ is indeed satisfied, with $ \varphi_{\log}^\kappa(t)$ given by \refeq{log-type-phi},
which is considered for $ t \le 0.9 $ and
$ \kappa = 1.8 $, and
\begin{align*}
w = (w_n), \quad \textup{with }
w_n =
\frac{4^\mykap}{\sqrt{n} (\log n)^{0.51}}, \quad n = 2, 3, \ldots \ .
\end{align*}
\begin{table}[h]
\hfill
\begin{tabular}{|| c | c |@{\hspace{5mm} } c | c ||}
\hline
\hline
$ \delta $
& $ 100 \myast \delta/\normb{f} $
& $ \ \normb{\upardeldel - \ust} $
& $ \ \normb{\upardeldel - \ust} \ / \varphi_{\log}^\kappa(\delta)  \ $
\\ \hline \hline
 $8.00 \myast 10^{-3}$ & $7.41 \myast 10^{-2}$ & $5.16 \myast 10^{-2}$ & $0.8786$ \\
 $4.00 \myast 10^{-3}$ & $3.71 \myast 10^{-2}$ & $4.31 \myast 10^{-2}$ & $0.9336$ \\
 $2.00 \myast 10^{-3}$ & $1.85 \myast 10^{-2}$ & $3.65 \myast 10^{-2}$ & $0.9776$ \\
 $1.00 \myast 10^{-3}$ & $9.27 \myast 10^{-3}$ & $3.12 \myast 10^{-2}$ & $1.0101$ \\
 $5.00 \myast 10^{-4}$ & $4.63 \myast 10^{-3}$ & $2.68 \myast 10^{-2}$ & $1.0328$ \\
 $2.50 \myast 10^{-4}$ & $2.32 \myast 10^{-3}$ & $2.32 \myast 10^{-2}$ & $1.0457$ \\
 $1.25 \myast 10^{-4}$ & $1.16 \myast 10^{-3}$ & $2.01 \myast 10^{-2}$ & $1.0483$ \\
 $6.25 \myast 10^{-5}$ & $5.79 \myast 10^{-4}$ & $1.75 \myast 10^{-2}$ & $1.0395$ \\
 $3.12 \myast 10^{-5}$ & $2.90 \myast 10^{-4}$ & $1.51 \myast 10^{-2}$ & $1.0193$ \\
 $1.56 \myast 10^{-5}$ & $1.45 \myast 10^{-4}$ & $1.30 \myast 10^{-2}$ & $0.9856$ \\
 $7.81 \myast 10^{-6}$ & $7.24 \myast 10^{-5}$ & $1.11 \myast 10^{-2}$ & $0.9361$ \\
 $3.91 \myast 10^{-6}$ & $3.62 \myast 10^{-5}$ & $9.26 \myast 10^{-3}$ & $0.8671$ \\
 $1.95 \myast 10^{-6}$ & $1.81 \myast 10^{-5}$ & $7.49 \myast 10^{-3}$ & $0.7728$ \\
\hline
\hline
\end{tabular}
\hfill \mbox{}
\caption{Numerical results for the a priori parameter choice strategy}
\label{tab:num1}
\end{table}

\subsection{Numerical results for the discrepancy principle}
\label{example2}
We next consider the discrepancy principle, \cf
Algorithm \ref{th:discrepancy-def}, with $ b = 4 $ and
for different values of $ \delta $.
It is in fact realized by the sequential version considered in Remark \ref{th:sdp},
with $ \theta = 10 $.
The numerical results are shown in Table \ref{tab:num2}.
The results presented in columns 3 and 5 confirm the statement of
Theorem \ref{th:aposteriori-convergence}.
The results in the last column are presented due to the discussion on the asymptotical behavior \refeq{infty}.
\begin{table}[h]
\hfill
\begin{tabular}{|| c | c | c | c | c | c ||}
\hline
\hline
$ \delta $
& $ 100 \myast \delta/\normb{\fst} $
& $ \ \normb{\upardeldel - \ust} $
& $ \ \pardel $
& $ \ \delta/\pardel^{1/4 } $
& $ \ \delta^2/\pardel $
\\ \hline \hline
 $1.00 \myast 10^{-3}$ & $9.27 \myast 10^{-3}$ & $4.09 \myast 10^{-2}$ & $1.00 \myast 10^{-5}$ & $1.78 \myast 10^{-2}$ & $  0.10$ \\
 $5.00 \myast 10^{-4}$ & $4.63 \myast 10^{-3}$ & $3.13 \myast 10^{-2}$ & $1.00 \myast 10^{-6}$ & $1.58 \myast 10^{-2}$ & $  0.25$ \\
 $2.50 \myast 10^{-4}$ & $2.32 \myast 10^{-3}$ & $2.44 \myast 10^{-2}$ & $1.00 \myast 10^{-7}$ & $1.41 \myast 10^{-2}$ & $  0.62$ \\
 $1.25 \myast 10^{-4}$ & $1.16 \myast 10^{-3}$ & $1.92 \myast 10^{-2}$ & $1.00 \myast 10^{-8}$ & $1.25 \myast 10^{-2}$ & $  1.56$ \\
 $6.25 \myast 10^{-5}$ & $5.79 \myast 10^{-4}$ & $1.51 \myast 10^{-2}$ & $1.00 \myast 10^{-9}$ & $1.11 \myast 10^{-2}$ & $  3.91$ \\
 $3.12 \myast 10^{-5}$ & $2.90 \myast 10^{-4}$ & $1.16 \myast 10^{-2}$ & $1.00 \myast 10^{-10}$ & $9.88 \myast 10^{-3}$ & $  9.77$ \\
 $1.56 \myast 10^{-5}$ & $1.45 \myast 10^{-4}$ & $1.17 \myast 10^{-2}$ & $1.00 \myast 10^{-10}$ & $4.94 \myast 10^{-3}$ & $  2.44$ \\
 $7.81 \myast 10^{-6}$ & $7.24 \myast 10^{-5}$ & $8.64 \myast 10^{-3}$ & $1.00 \myast 10^{-11}$ & $4.39 \myast 10^{-3}$ & $  6.10$ \\
 $3.91 \myast 10^{-6}$ & $3.62 \myast 10^{-5}$ & $5.62 \myast 10^{-3}$ & $1.00 \myast 10^{-12}$ & $3.91 \myast 10^{-3}$ & $ 15.26$ \\
 $1.95 \myast 10^{-6}$ & $1.81 \myast 10^{-5}$ & $2.48 \myast 10^{-3}$ & $1.00 \myast 10^{-13}$ & $3.47 \myast 10^{-3}$ & $ 38.15$ \\
\hline
\hline
\end{tabular}
\hfill \mbox{}
\caption{Numerical results for the discrepancy principle}
\label{tab:num2}
\end{table}

\section*{Acknowledgment}
This research has been supported by German Research Foundation
(DFG grant HO 1454/12-1) under the auspices of the joint Austrian--German project
``Novel Error Measures and Source Conditions of Regularization Methods for
Inverse Problems (SCIP)'' with the University of Vienna
(PI: Prof.~Dr.~Otmar Scherzer) according to D-A-CH Lead Agency
Agreement.

The authors are grateful for suggestions of two referees,
which led to an improved presentation.

\bibliography{HP}
\end{document}